\documentclass[11pt]{amsart}
\usepackage{amsmath,amssymb}
\newtheorem{theorem}{Theorem}[section]
\newtheorem{corollary}[theorem]{Corollary}
\newtheorem{proposition}[theorem]{Proposition} 
\newtheorem{lemma}[theorem]{Lemma}

\begin{document}

\title[Minimal Lagrangian Diffeomorphisms]{Minimal Lagrangian Diffeomorphisms between Domains in the Hyperbolic Plane}
\author{Simon Brendle}

\maketitle

\section{Introduction}

This paper is concerned with the boundary regularity of minimal Lagrangian diffeomorphisms. The notion of a minimal Lagrangian diffeomorphism was introduced by R.~Schoen. In \cite{Schoen}, Schoen proved an existence and uniqueness result for minimal Lagrangian diffeomorphisms between hyperbolic surfaces:

\begin{theorem}[R.~Schoen \cite{Schoen}] 
\label{thm.1}
Let $N$ be a compact surface of genus greater than $1$, and let $g,\tilde{g}$ be a pair of hyperbolic metrics on $N$. Then there exists a unique diffeomorphism $f: N \to N$ with the following properties: 
\begin{itemize}
\item[(i)] $f$ is area-preserving 
\item[(ii)] $f$ is homotopic to the identity 
\item[(iii)] The graph of $f$ is a minimal submanifold of $(N,g) \times (N,\tilde{g})$
\end{itemize}
\end{theorem}

Theorem \ref{thm.1} was subsequently generalized by Y.I.~Lee \cite{Lee}. M.T.~Wang \cite{Wang1} gave an alternative proof of the existence part of Theorem \ref{thm.1} using mean curvature flow.

In this paper, we study an analogous problem for surfaces with boundary. Throughout this paper, we will assume that $N$ is a complete, simply-connected surface of constant curvature $\kappa \leq 0$. Suppose that $\Omega$ and $\tilde{\Omega}$ are domains in $N$ with smooth boundary, and let $f$ be a diffeomorphism from $\Omega$ to $\tilde{\Omega}$. We will say that $f$ is a minimal Lagrangian diffeomorphism if the following conditions are satisfied: 
\begin{itemize}
\item[(i)] $f$ is area-preserving 
\item[(ii)] $f$ is orientation-preserving 
\item[(iii)] The graph of $f$ is a minimal submanifold of $N \times N$
\end{itemize} 
The case $\kappa = 0$ is somewhat special. In this case, the existence of a minimal Lagrangian diffeomorphism from $\Omega$ to $\tilde{\Omega}$ is closely related to the solvability of the second boundary value problem for the Monge-Amp\`ere equation (cf. \cite{Wolfson}). To describe the connection between the two problems, we consider two domains $\Omega,\tilde{\Omega} \subset \mathbb{R}^2$ and a minimal Lagrangian diffeomorphism $f: \Omega \to \tilde{\Omega}$. By definition, the graph of $f$ is a minimal Lagrangian submanifold of $\mathbb{R}^2 \times \mathbb{R}^2$. Consequently, the graph of $f$ has constant Lagrangian angle. This implies that $f$ is the composition of a gradient mapping $x \mapsto \nabla u(x)$ and a rotation. Since $f$ is area-preserving and orientation-preserving, the function $u: \Omega \to \mathbb{R}$ is a solution of the Monge-Amp\`ere equation $\det D^2 u(x) = 1$.

P.~Delano\"e \cite{Delanoe} has obtained an existence result for the second boundary value problem for the Monge-Amp\`ere equation in dimension $2$. This result was extended to higher dimensions by L.~Caffarelli \cite{Caffarelli} and J.~Urbas \cite{Urbas1}. The following result is an immediate consequence of Delano\"e's existence theorem:

\begin{theorem}[P.~Delano\"e \cite{Delanoe}] 
\label{thm.2}
Let $\Omega$ and $\tilde{\Omega}$ be strictly convex domains in $\mathbb{R}^2$ with smooth boundary. Assume that $\Omega$ and $\tilde{\Omega}$ have the same area. Then there exists a minimal Lagrangian diffeomorphism from $\Omega$ to $\tilde{\Omega}$.
\end{theorem}

We point out that the convexity of both domains $\Omega,\tilde{\Omega}$ is essential. J.~Urbas \cite{Urbas2} has recently constructed an example of a non-convex domain of area $\pi$ that does not admit a minimal Lagrangian diffeomorphism to the unit disk.

We now return to the general case ($\kappa \leq 0$). Note that the link between minimal Lagrangian diffeomorphisms and solutions of the Monge-Amp\`ere equation breaks down in this setting. Nonetheless, we have the following existence and uniqueness result:

\begin{theorem} 
\label{existence.uniqueness}
Let $\Omega$ and $\tilde{\Omega}$ be strictly convex domains in $N$ with smooth boundary. Assume that $\Omega$ and $\tilde{\Omega}$ have the same area. Given any point $\overline{p} \in \partial \Omega$ and any point $\overline{q} \in \partial \tilde{\Omega}$, there exists a unique minimal Lagrangian diffeomorphism from $\Omega$ to $\tilde{\Omega}$ that maps $\overline{p}$ to $\overline{q}$.
\end{theorem}

In order to prove Theorem \ref{existence.uniqueness}, we deform $\Omega$ and $\tilde{\Omega}$ to the flat unit disk $\mathbb{B}^2 \subset \mathbb{R}^2$, and apply the continuity method. This requires a-priori estimates for minimal Lagrangian diffeomorphisms from $\Omega$ to $\tilde{\Omega}$.

\begin{theorem}
\label{a.priori.estimates}
Let $\Omega$ and $\tilde{\Omega}$ be strictly convex domains in $N$ with smooth boundary. Suppose that $f: \Omega \to \tilde{\Omega}$ is a minimal Lagrangian diffeomorphism. Then $f$ is uniformly bounded in $C^m$; more precisely, we have $\|f\|_{C^m} \leq C$, where $C = C(m,\Omega,\tilde{\Omega})$ is a constant that depends only on $m$ and the domains $\Omega$ and $\tilde{\Omega}$.
\end{theorem}

The proof of Theorem \ref{a.priori.estimates} will occupy Sections 2 -- 6. In Section 2, we construct boundary defining functions for $\Omega$ and $\tilde{\Omega}$ which are uniformly convex. Moreover, we establish some basic estimates involving the boundary defining functions. In Section 3, we introduce tools from complex geometry. In Section 4, we use these ideas to estimate the singular values of $Df_p$ for all boundary points $p \in \partial \Omega$. In Section 5, we employ an argument due to M.T.~Wang to obtain uniform bounds for the singular values of $Df_p$ for all $p \in \Omega$. In Section 6, we show that $f$ is bounded in $C^{1,\alpha}$. In Section 7, we show that the linearized operator is invertible. This precludes bifurcations. Finally, in Section 8, we show that every minimal Lagrangian diffeomorphism from the flat unit disk to itself is a rotation. This follows from a uniqueness result, due to P.~Delano\"e \cite{Delanoe}, for the second boundary value problem for the Monge-Amp\`ere equation.

The author is grateful to Professor Richard Schoen and Professor Leon Simon for discussions. This project was supported by the Alfred P. Sloan Foundation and by the National Science Foundation under grant DMS-0605223.

\section{The boundary defining functions}

As above, we assume that $\Omega$ and $\tilde{\Omega}$ are strictly convex domains in $N$ with smooth boundary. We begin by constructing a boundary defining function for the domain $\Omega$ which is uniformly convex:

\begin{proposition}
\label{boundary.defining.functions}
There exists a smooth function $h: \Omega \to \mathbb{R}$ with the following properties: 
\begin{itemize}
\item $h$ is uniformly convex
\item For each point $p \in \partial \Omega$, we have $h(p) = 0$ and $|\nabla h(p)| = 1$
\item If $s$ is sufficiently close to $\inf_\Omega h$, then the sub-level set $\{p \in \Omega: h(p) \leq s\}$ is a geodesic disk
\end{itemize}
Similarly, we can find a smooth function $\tilde{h}: \tilde{\Omega} \to \mathbb{R}$ such that: 
\begin{itemize}
\item $\tilde{h}$ is uniformly convex
\item For each point $q \in \partial \tilde{\Omega}$, we have $\tilde{h}(q) = 0$ and $|\nabla \tilde{h}_q| = 1$
\item If $s$ is sufficiently close to $\inf_{\tilde{\Omega}} \tilde{h}$, then the sub-level set $\{q \in \tilde{\Omega}: \tilde{h}(q) \leq s\}$ is a geodesic disk
\end{itemize}
\end{proposition}

\textbf{Proof.} 
We will only prove the assertion for the domain $\Omega$. Let $p_0$ be an arbitrary point in the interior of $\Omega$. We define a function $h_1$ by 
\[h_1(p) = \frac{d(p,\partial \Omega)^2}{4 \, \text{\rm diam}(\Omega)} - d(p,\partial \Omega).\] 
Since $\Omega$ is strictly convex, there exists a positive real number $\varepsilon$ such that $h_1$ is smooth and uniformly convex for $d(p,\partial \Omega) < \varepsilon$. We assume that $\varepsilon$ is chosen so that $d(p_0,\partial \Omega) > \varepsilon$. We next define a function $h_2$ by 
\[h_2(p) = \frac{\varepsilon \, d(p_0,p)^2}{4 \, \text{\rm diam}(\Omega)^2} - \frac{\varepsilon}{2}.\] 
The function $h_2$ is smooth and uniformly convex by the Hessian comparison theorem. For each point $p \in \partial \Omega$, we have $h_1(p) = 0$ and $h_2(p) \leq -\frac{\varepsilon}{4}$. Moreover, for $d(p,\partial \Omega) \geq \varepsilon$, we have $h_1(p) \leq -\frac{3\varepsilon}{4}$ and $h_2(p) \geq -\frac{\varepsilon}{2}$. 

We now define 
\[h(p) = \frac{h_1(p) + h_2(p)}{2} + \Phi \Big ( \frac{h_1(p) - h_2(p)}{2} \Big ),\] 
where $\Phi: \mathbb{R} \to \mathbb{R}$ is a smooth function satisfying $\Phi''(s) \geq 0$ for all $s \in \mathbb{R}$ and $\Phi(s) = |s|$ for $|s| \geq \frac{\varepsilon}{16}$. It is easy to see that $h$ is smooth and uniformly convex for $d(p,\partial \Omega) < \varepsilon$. Since $h$ agrees with $h_2$ for $d(p,\partial \Omega) \geq \varepsilon$, we conclude that $h$ is smooth and uniformly convex in all of $\Omega$. Moreover, $h$ agrees with $h_1$ in a neighborhood of $\partial \Omega$. 
Thus, we conclude that $h(p) = 0$ and $|\nabla h_p| = 1$ for all $p \in \partial \Omega$. 

It remains to verify the last statement. It is easy to see that $h(p) \geq h_2(p) \geq -\frac{\varepsilon}{2}$ for all $p \in \Omega$. Since $h(p_0) = h_2(p_0) = -\frac{\varepsilon}{2}$, it follows that $\inf_\Omega h = -\frac{\varepsilon}{2}$. Moreover, if $s$ is a real number satisfying \[-\frac{\varepsilon}{2} < s < \frac{\varepsilon \, (d(p_0,\partial \Omega) - \varepsilon)^2}{4 \, \text{\rm diam}(\Omega)^2} - \frac{\varepsilon}{2},\] then the set $\{p \in \Omega: h(p) \leq s\}$ is a geodesic disk. This completes the proof of Proposition \ref{boundary.defining.functions}. \\

Since $h$ and $\tilde{h}$ are uniformly convex, we can find a positive constant $\theta$ such that
\[\theta \, |w|^2 \leq (\text{\rm Hess} \, h)_p(w,w) \leq \frac{1}{\theta} \, |w|^2\]
for all $p \in \Omega$ and $w \in T_p N$ and 
\[\theta \, |\tilde{w}|^2 \leq (\text{\rm Hess} \, \tilde{h})_q(\tilde{w},\tilde{w}) \leq \frac{1}{\theta} \, |\tilde{w}|^2\] 
for all $q \in \tilde{\Omega}$ and $\tilde{w} \in T_q N$. 

Suppose now that $f: \Omega \to \tilde{\Omega}$ is a minimal Lagrangian diffeomorphism. Let 
\[\Sigma = \{(p,f(p)): p \in \Omega\}\] 
be the graph of $f$. By definition, $\Sigma$ is a minimal submanifold of the product manifold $M = N \times N$. We define two functions $H,\tilde{H}: \Sigma \to \mathbb{R}$ by $H(p,f(p)) = h(p)$ and $\tilde{H}(p,f(p)) = \tilde{h}(f(p))$. 

\begin{proposition} 
\label{laplacian.of.H}
The function $H$ satisfies $\theta \leq \Delta_\Sigma H \leq \frac{1}{\theta}$.
Similarly, the function $\tilde{H}$ satisfies $\theta \leq \Delta_\Sigma \tilde{H} \leq \frac{1}{\theta}$.
\end{proposition}

\textbf{Proof.} Fix a point $(p,f(p)) \in \Sigma$. We can find an orthonormal basis $\{v_1,v_2\}$ of $T_p N$ such that $\big [ Df_p^* \, Df_p \big ] \, v_k = \lambda_k^2 \, v_k$, where $\lambda_1,\lambda_2$ are positive real numbers 
satisfying $\lambda_1 \lambda_2 = 1$. Since $\Sigma$ is a minimal submanifold of $M$, the Laplacian of $H$ at $(p,f(p))$ is given by 
\begin{align*} 
\Delta_\Sigma H 
&= \text{\rm tr} \Big[ (I + Df_p^* \, Df_p)^{-1} \, (\text{\rm Hess} \, h)_p \Big ] \\ 
&= \sum_{k=1}^2 \frac{1}{1 + \lambda_k^2} \, (\text{\rm Hess} \, h)_p(v_k,v_k). 
\end{align*} 
By assumption, we have $\theta \leq (\text{\rm Hess} \, h)_p(v_k,v_k) \leq \frac{1}{\theta}$ 
for $k = 1,2$. Moreover, the relation $\lambda_1 \lambda_2 = 1$ implies \[\frac{1}{1 + \lambda_1^2} + \frac{1}{1 + \lambda_2^2} = 1.\] Thus, we conclude that $\theta \leq \Delta_\Sigma H \leq \frac{1}{\theta}$. The inequality $\theta \leq \Delta_\Sigma \tilde{H} \leq \frac{1}{\theta}$ follows from an analogous argument. \\

\begin{proposition} 
\label{estimate.for.the.boundary.defining.functions.1}
We have $-\theta^2 \, h(p) \leq -\tilde{h}(f(p)) \leq -\frac{1}{\theta^2} \, h(p)$ for all $p \in \Omega$. 
\end{proposition}

\textbf{Proof.} It follows from Proposition \ref{laplacian.of.H} that the functions
$\theta^2 \, H - \tilde{H}$ and $\theta^2 \, \tilde{H} - H$ are superharmonic. Since both $H$ and $\tilde{H}$ vanish along the boundary of $\Sigma$, we conclude that $-\theta^2 \, H \leq -\tilde{H} \leq -\frac{1}{\theta^2} \, H$ by the maximum principle. From this, the assertion follows. \\

\begin{corollary} 
\label{estimate.for.the.boundary.defining.functions.2}
We have $\theta^2 \leq \langle Df_p(\nabla h_p),\nabla \tilde{h}_{f(p)} \rangle \leq \frac{1}{\theta^2}$ for all $p \in \partial \Omega$. 
\end{corollary}

\section{Tools from complex geometry}

As in the previous section, we assume that $f: \Omega \to \tilde{\Omega}$ is a minimal Lagrangian diffeomorphism. Fix a complex structure $J$ on $N$. We define a complex structure on the product $M = N \times N$ by $J_{(p,q)}(w,\tilde{w}) = (J_p w,-J_q \tilde{w})$ for all vectors $w \in T_p N$ and $\tilde{w} \in T_q N$. Since $f$ is area-preserving and orientation-preserving, the graph $\Sigma = \{(p,f(p)): p \in \Omega\}$ is a Lagrangian submanifold of $M$.

For each point $p \in \Omega$, we define a linear isometry $Q_p: T_p N \to T_{f(p)} N$ by
\[Q_p = Df_p \, \big [ Df_p^* \, Df_p \big ]^{-\frac{1}{2}}.\] 
It is easy to see that $Q_p: T_p N \to T_{f(p)} N$ is orientation-preserving. This implies $J_{f(p)} \, Q_p = Q_p \, J_p$ for all $p \in \Omega$. 

For each point $p \in \Omega$, we define a bilinear form $\sigma: T_{(p,f(p))} M \times T_{(p,f(p))} M \to \mathbb{C}$ by 
\begin{align*} 
\sigma \big ( (w_1,\tilde{w}_1),(w_2,\tilde{w}_2) \big ) 
&= i \, \langle Q_p(w_1),\tilde{w}_2 \rangle + \langle Q_p(J_p w_1),\tilde{w}_2 \rangle \\ 
&- i \, \langle Q_p(w_2),\tilde{w}_1 \rangle - \langle Q_p(J_p w_2),\tilde{w}_1 \rangle 
\end{align*} 
for $w_1,w_2 \in T_p N$ and $\tilde{w}_1,\tilde{w}_2 \in T_{f(p)} N$. \\

\begin{lemma}
\label{complex.volume.form}
We have $\sigma(W_2,W_1) = -\sigma(W_1,W_2)$ and $\sigma(JW_1,W_2) = i \, \sigma(W_1,W_2)$ for all $W_1,W_2 \in T_{(p,f(p))} M$.
\end{lemma}

\textbf{Proof.} 
The first property is trivial. To prove the second property, we observe that 
\begin{align*} 
\sigma \big ( (J_p w_1,-J_{f(p)} \tilde{w}_1),(w_2,\tilde{w}_2) \big ) 
&= i \, \langle Q_p(J_p w_1),\tilde{w}_2 \rangle - \langle Q_p(w_1),\tilde{w}_2 \rangle \\ 
&+ i \, \langle Q_p(w_2),J_{f(p)} \tilde{w}_1 \rangle + \langle Q_p(J_p w_2),J_{f(p)} \tilde{w}_1 \rangle \\ 
&= i \, \langle Q_p(J_p w_1),\tilde{w}_2 \rangle - \langle Q_p(w_1),\tilde{w}_2 \rangle \\ 
&- i \, \langle Q_p(J_p w_2),\tilde{w}_1 \rangle + \langle Q_p(w_2),\tilde{w}_1 \rangle \\ 
&= i \, \sigma \big ( (w_1,\tilde{w}_1),(w_2,\tilde{w}_2) \big ) 
\end{align*} 
for all vectors $w_1,w_2 \in T_p N$ and $\tilde{w}_1,\tilde{w}_2 \in T_{f(p)} N$. \\

\begin{lemma}
\label{lagrangian.angle}
If $\{e_1,e_2\}$ is an orthonormal basis of $T_{(p,f(p))} \Sigma$, then $\sigma(e_1,e_2) = \pm 1$. 
\end{lemma}

\textbf{Proof.} 
Since $\sigma$ is anti-symmetric, it is enough to prove the assertion for one particular orthonormal basis of $T_{(p,f(p))} \Sigma$. To that end, we choose an orthonormal basis $\{v_1,v_2\}$ of $T_p N$ such that $\big [ Df_p^* \, Df_p \big ] \, v_k = \lambda_k^2 \, v_k$, where $\lambda_1,\lambda_2$ are positive real numbers satisfying $\lambda_1 \lambda_2 = 1$. Since $Q_p$ is an isometry, we have 
\[\langle Q_p(J_p v_1),Q_p(v_2) \rangle = -\langle Q_p(J_p v_2),Q_p(v_1) \rangle = \pm 1\] 
and 
\[\langle Q_p(v_1),Q_p(v_2) \rangle = 0.\] 
Moreover, the relation $\big [ Df_p^* \, Df_p \big ] \, v_k = \lambda_k^2 \, v_k$ implies $Df_p(v_k) = \lambda_k \, Q_p(v_k)$. We now define 
\[e_k = \frac{1}{\sqrt{1 + \lambda_k^2}} \, (v_k,Df_p(v_k)) = \frac{1}{\sqrt{1 + \lambda_k^2}} \, (v_k,\lambda_k \, Q_p(v_k))\] 
for $k = 1,2$. Clearly, $\{e_1,e_2\}$ is an orthonormal basis of $T_{(p,f(p))} \Sigma$. By definition of $\sigma$, we have 
\begin{align*} 
&\sigma(e_1,e_2) \\ 
&= \frac{1}{\sqrt{1 + \lambda_1^2}} \, \frac{1}{\sqrt{1 + \lambda_2^2}} \, \Big [ i\lambda_2 \, \langle Q_p(v_1),Q_p(v_2) \rangle + \lambda_2 \, \langle Q_p(J_p v_1),Q_p(v_2) \rangle \\ &\hspace{38mm} - i\lambda_1 \, \langle Q_p(v_2),Q_p(v_1) \rangle - \lambda_1 \, \langle Q_p(J_p v_2),Q_p(v_1) \rangle \Big ] \\ 
&= \pm \frac{1}{\sqrt{1 + \lambda_1^2}} \, \frac{1}{\sqrt{1 + \lambda_2^2}} \, (\lambda_1 + \lambda_2) \\ 
&= \pm 1. 
\end{align*}
This proves the assertion. \\

We next show that $\sigma$ is parallel with respect to the Levi-Civita connection on $M$. To fix notation, we denote by $TM|_\Sigma$ the restriction of the tangent bundle $TM$ to $\Sigma$.

\begin{proposition}
\label{gradient.of.psi}
Let $W_1,W_2$ be sections of the vector bundle $TM|_\Sigma$. We define a function $\psi: \Sigma \to \mathbb{C}$ by $\psi = \sigma(W_1,W_2)$. Then \[V(\psi) = \sigma(\nabla_V^M W_1,W_2) + \sigma(W_1,\nabla_V^M W_2)\] for all $V \in T\Sigma$.
\end{proposition}

\textbf{Proof.} Fix a tangent vector field $V$ along $\Sigma$, and let 
\[\tau(W_1,W_2) := \sigma(\nabla_V^M W_1,W_2) + \sigma(W_1,\nabla_V^M W_2) - V(\sigma(W_1,W_2)).\] 
It is easy to see that $\tau(W_1,W_2)$ is a tensor. It follows from Lemma \ref{complex.volume.form} that 
$\tau(W_2,W_1) = -\tau(W_1,W_2)$ and $\tau(JW_1,W_2) = i \, \tau(W_1,W_2)$. Hence, it suffices to show that $\tau(e_1,e_2) = 0$, where $\{e_1,e_2\}$ is a local orthonormal frame on $\Sigma$. By Lemma \ref{lagrangian.angle}, we have $\sigma(e_1,e_2) = \pm 1$. This implies $V(\sigma(e_1,e_2)) = 0$. Moreover, we have $\sigma(\nabla_V^\Sigma e_1,e_2) = 0$ since $\nabla_V^\Sigma e_1$ is a multiple of $e_2$. Similarly, $\sigma(e_1,\nabla_V^\Sigma e_2) = 0$ since $\nabla_V^\Sigma e_2$ is a multiple of $e_1$. Putting these facts together, we obtain 
\begin{align*}
\tau(e_1,e_2) 
&= \sigma(\nabla_V^M e_1 - \nabla_V^\Sigma e_1,e_2) + \sigma(e_1,\nabla_V^M e_2 - \nabla_V^\Sigma e_2) \\ 
&= \sigma(I\!I(e_1,V),e_2) + \sigma(e_1,I\!I(e_2,V)) \\ 
&= \sigma(Je_1,e_2) \, \langle I\!I(e_1,V),Je_1 \rangle + \sigma(e_1,Je_1) \, \langle I\!I(e_2,V),Je_1 \rangle \\ 
&+ \sigma(Je_2,e_2) \, \langle I\!I(e_1,V),Je_2 \rangle + \sigma(e_1,Je_2) \, \langle I\!I(e_2,V),Je_2 \rangle \\ 
&= i \, \sigma(e_1,e_2) \, (\langle I\!I(e_1,V),Je_1 \rangle + \langle I\!I(e_2,V),Je_2 \rangle) \\ 
&= i \, \sigma(e_1,e_2) \, (\langle I\!I(e_1,e_1),JV \rangle + \langle I\!I(e_2,e_2),JV \rangle), 
\end{align*}
where $I\!I$ denotes the second fundamental form of $\Sigma$. Since \[I\!I(e_1,e_1) + I\!I(e_2,e_2) = 0,\] 
we conclude that $\tau(e_1,e_2) = 0$. This completes the proof of Proposition \ref{gradient.of.psi}. \\

Given a smooth vector field $W$ on $M$, we denote by 
\[\nabla_{V_1,V_2}^{M,2} W = \nabla_{V_1}^M \nabla_{V_2}^M W - \nabla_{\nabla_{V_1}^M V_2}^M W\] 
the second order covariant derivative of a vector field $W$ with respect to the Levi-Civita connection on $M$. 

We next compute the Hessian of a function of the form $\psi = \sigma(W_1,W_2)$, where $W_1,W_2$ are smooth vector fields on $M$.

\begin{proposition}
\label{hessian.of.psi}
Let $W_1,W_2$ be smooth vector fields on $M$. We define a function $\psi: \Sigma \to \mathbb{C}$ by $\psi = \sigma(W_1,W_2)$. Then 
\begin{align*} 
(\text{\rm Hess}_\Sigma \, \psi)(V_1,V_2) 
&= \sigma(\nabla_{V_1,V_2}^{M,2} W_1,W_2) + \sigma(W_1,\nabla_{V_1,V_2}^{M,2} W_2) \\
&+ \sigma(\nabla_{V_1}^M W_1,\nabla_{V_2}^M W_2) + \sigma(\nabla_{V_2}^M W_1,\nabla_{V_1}^M W_2) \\ 
&+ \sigma(\nabla_{I\!I(V_1,V_2)}^M W_1,W_2) + \sigma(W_1,\nabla_{I\!I(V_1,V_2)}^M W_2) 
\end{align*} 
for all $V_1,V_2 \in T\Sigma$.
\end{proposition}

\textbf{Proof.} 
Suppose that $V_1,V_2$ are tangent vector fields along $\Sigma$. It follows from the previous proposition that 
\[V_2(\psi) = \sigma(\nabla_{V_2}^M W_1,W_2) + \sigma(W_1,\nabla_{V_2}^M W_2).\] 
This implies 
\begin{align*} 
V_1(V_2(\psi)) 
&= \sigma(\nabla_{V_1}^M \nabla_{V_2}^M W_1,W_2) + \sigma(W_1,\nabla_{V_1}^M \nabla_{V_2}^M W_2) \\ 
&+ \sigma(\nabla_{V_1}^M W_1,\nabla_{V_2}^M W_2) + \sigma(\nabla_{V_2}^M W_1,\nabla_{V_1}^M W_2). 
\end{align*}
Thus, we conclude that 
\begin{align*} 
(\text{\rm Hess}_\Sigma \, \psi)(V_1,V_2) 
&= V_1(V_2(\psi)) - (\nabla_{V_1}^\Sigma V_2)(\psi) \\ 
&= \sigma(\nabla_{V_1}^M \nabla_{V_2}^M W_1,W_2) + \sigma(W_1,\nabla_{V_1}^M \nabla_{V_2}^M W_2) \\ 
&+ \sigma(\nabla_{V_1}^M W_1,\nabla_{V_2}^M W_2) + \sigma(\nabla_{V_2}^M W_1,\nabla_{V_1}^M W_2) \\ 
&- \sigma(\nabla_{\nabla_{V_1}^\Sigma V_2}^M W_1,W_2) - \sigma(W_1,\nabla_{\nabla_{V_1}^\Sigma V_2}^M W_2). 
\end{align*}
Using the identity $\nabla_{V_1}^M V_2 - \nabla_{V_1}^\Sigma V_2 = I\!I(V_1,V_2)$, we obtain 
\begin{align*} 
(\text{\rm Hess}_\Sigma \, \psi)(V_1,V_2) 
&= \sigma(\nabla_{V_1}^M \nabla_{V_2}^M W_1,W_2) + \sigma(W_1,\nabla_{V_1}^M \nabla_{V_2}^M W_2) \\ 
&+ \sigma(\nabla_{V_1}^M W_1,\nabla_{V_2}^M W_2) + \sigma(\nabla_{V_2}^M W_1,\nabla_{V_1}^M W_2) \\ 
&- \sigma(\nabla_{\nabla_{V_1}^M V_2}^M W_1,W_2) - \sigma(W_1,\nabla_{\nabla_{V_1}^M V_2}^M W_2) \\ 
&+ \sigma(\nabla_{I\!I(V_1,V_2)}^M W_1,W_2) + \sigma(W_1,\nabla_{I\!I(V_1,V_2)}^M W_2). 
\end{align*}
From this, the assertion follows. \\

\begin{corollary}
\label{laplacian.of.psi}
Let $W_1,W_2$ be smooth vector fields on $M$. As above, we define a function $\psi: \Sigma \to \mathbb{C}$ by $\psi = \sigma(W_1,W_2)$. Then 
\begin{align*} 
\Delta_\Sigma \psi 
&= \sum_{k=1}^2 \sigma(\nabla_{e_k,e_k}^{M,2} W_1,W_2) + \sum_{k=1}^2 \sigma(W_1,\nabla_{e_k,e_k}^{M,2} W_2) \\ 
&+ 2 \, \sum_{k=1}^2 \sigma(\nabla_{e_k}^M W_1,\nabla_{e_k}^M W_2). 
\end{align*}
\end{corollary}

\textbf{Proof.} 
This follows immediately from Proposition \ref{hessian.of.psi} and the relation $\sum_{k=1}^2 I\!I(e_k,e_k) = 0$. \\

\section{The boundary gradient estimate}

We define a vector field $\xi$ on $\Omega$ by $\xi = \nabla h$; similarly, we define a vector field $\tilde{\xi}$ on $\tilde{\Omega}$ by $\tilde{\xi} = \nabla \tilde{h}$. We next define a function $\varphi: \Sigma \to \mathbb{R}$ by
\[\varphi(p,f(p)) = \langle Q_p(\xi_p),\tilde{\xi}_{f(p)} \rangle\] 
for all $p \in \Omega$.

\begin{proposition} 
\label{properties.of.varphi}
The gradient of the function $\varphi: \Sigma \to \mathbb{R}$ is given by 
\begin{align*} 
\big \langle \nabla^\Sigma \varphi_{(p,f(p))},(v,Df_p(v)) \big \rangle 
&= \langle Q_p(\nabla_v \xi),\tilde{\xi}_{f(p)} \rangle \\ 
&+ \langle Q_p(\xi_p),\nabla_{Df_p(v)} \tilde{\xi} \rangle 
\end{align*} 
for all $p \in \Omega$ and $v \in T_p N$. Moreover, there exists a constant $C_1$, depending only on $h$ and $\tilde{h}$, such that 
$|\Delta_\Sigma \varphi| \leq C_1$ at each point $(p,f(p)) \in \Sigma$. 
\end{proposition}

\textbf{Proof.} 
We define two vector fields $W_1$ and $W_2$ on $\Omega \times \tilde{\Omega} \subset M$ by 
\[(W_1)_{(p,q)} = (\xi_p,0) \in T_p N \times T_q N\] 
and 
\[(W_2)_{(p,q)} = (0,\tilde{\xi}_q) \in T_p N \times T_q N\] 
for all points $(p,q) \in \Omega \times \tilde{\Omega}$. As in the previous section, we define a function $\psi: \Sigma \to \mathbb{C}$ by $\psi = \sigma(W_1,W_2)$. This implies 
\[\psi(p,f(p)) = i \, \langle Q_p(\xi_p),\tilde{\xi}_{f(p)} \rangle + \langle Q_p(J_p \xi_p),\tilde{\xi}_{f(p)} \rangle\] 
for all $p \in \Omega$. Hence, the function $\varphi$ is the imaginary part of $\psi$. Using Proposition \ref{gradient.of.psi}, we obtain 
\begin{align*} 
\big \langle \nabla^\Sigma \psi_{(p,f(p))},(v,Df_p(v)) \big \rangle 
&= i \, \langle Q_p(\nabla_v \xi),\tilde{\xi}_{f(p)} \rangle + \langle Q_p(J_p \nabla_v \xi),\tilde{\xi}_{f(p)} \rangle \\ 
&+ i \, \langle Q_p(\xi_p),\nabla_{Df_p(v)} \tilde{\xi} \rangle + \langle Q_p(J_p \xi_p),\nabla_{Df_p(v)} \tilde{\xi} \rangle 
\end{align*} 
for all $p \in \Omega$ and $v \in T_p N$. Since $\varphi = \text{\rm Im}(\psi)$, the first statement follows. Moreover, it follows from Corollary \ref{laplacian.of.psi} that $|\Delta_\Sigma \psi| \leq C_1$ for some constant $C_1$. Since $\varphi = \text{\rm Im}(\psi)$, we conclude that $|\Delta_\Sigma \varphi| \leq C_1$. \\

\begin{proposition} 
\label{varphi.is.positive} 
We have $\varphi(p,f(p)) > 0$ for all $p \in \partial \Omega$.
\end{proposition}

\textbf{Proof.}
Fix a point $p \in \partial \Omega$. By definition of $Q_p$, we have 
\begin{align*} 
\big \langle Q_p(Df_p^*(\tilde{\xi}_{f(p)})),\tilde{\xi}_{f(p)} \big \rangle &= \big \langle Df_p^*(\tilde{\xi}_{f(p)}),Q_p^*(\tilde{\xi}_{f(p)}) \big \rangle \\ 
&= \big \langle [Df_p^* \, Df_p]^{\frac{1}{2}} \, Q_p^*(\tilde{\xi}_{f(p)}),Q_p^*(\tilde{\xi}_{f(p)}) \big \rangle \\ 
&> 0. 
\end{align*} 
On the other hand, the vector $Df_p^*(\tilde{\xi}_{f(p)})$ is a positive multiple of $\xi_p$. Thus, we conclude that $\langle Q_p(\xi_p),\tilde{\xi}_{f(p)} \rangle > 0$, as claimed. \\

\begin{proposition}
\label{lower.bound.for.varphi}
We have $\varphi(p,f(p)) \geq \frac{\theta^2}{C_1}$ for all $p \in \partial \Omega$.
\end{proposition}

\textbf{Proof.} 
It follows from Proposition \ref{properties.of.varphi} that the function $\varphi - \frac{C_1}{\theta} \, H$ is superharmonic. Hence, there exists a point $p_0 \in \partial \Omega$ such that 
\[\inf_{p \in \Omega} \Big ( \varphi(p,f(p)) - \frac{C_1}{\theta} \, H(p,f(p)) \Big ) 
= \inf_{p \in \partial \Omega} \varphi(p,f(p)) = \varphi(p_0,f(p_0)).\] 
By the Hopf boundary point lemma, there exists a real number $\mu \geq 0$ such that 
\[\nabla^\Sigma \varphi = \Big ( \frac{C_1}{\theta} - \mu \Big ) \, \nabla^\Sigma H\] 
at $(p_0,f(p_0))$. By Proposition \ref{properties.of.varphi}, we have 
\begin{align*} 
\big \langle \nabla^\Sigma \varphi,(v,Df_p(v)) \big \rangle 
&= \langle Q_p(\nabla_v \xi),\tilde{\xi}_{f(p)} \rangle \\ 
&+ \langle Q_p(\xi_p),\nabla_{Df_p(v)} \tilde{\xi} \rangle \\ 
&= (\text{\rm Hess} \, h)_p \big ( v,Q_p^*(\tilde{\xi}_{f(p)}) \big ) \\ 
&+ (\text{\rm Hess} \, \tilde{h})_{f(p)} \big ( Q_p(\xi_p),Df_p(v) \big ) 
\end{align*} 
for all $p \in \Omega$ and all $v \in T_p N$. This implies 
\begin{align*} 
\Big ( \frac{C_1}{\theta} - \mu \Big ) \, \langle \xi_{p_0},v \rangle 
&= (\text{\rm Hess} \, h)_{p_0} \big ( v,Q_{p_0}^*(\tilde{\xi}_{f(p_0)}) \big ) \\ 
&+ (\text{\rm Hess} \, \tilde{h})_{f(p_0)} \big ( Q_{p_0}(\xi_{p_0}),Df_{p_0}(v) \big ) 
\end{align*} 
for all $v \in T_{p_0} N$. Hence, if we put $v = Q_{p_0}^*(\tilde{\xi}_{f(p_0)})$, then we obtain 
\begin{align*} 
\Big ( \frac{C_1}{\theta} - \mu \Big ) \, \varphi(p_0,f(p_0)) 
&= (\text{\rm Hess} \, h)_{p_0} \big ( Q_{p_0}^*(\tilde{\xi}_{f(p_0)}),Q_{p_0}^*(\tilde{\xi}_{f(p_0)}) \big ) \\ 
&+ (\text{\rm Hess} \, \tilde{h})_{f(p_0)} \big ( Q_{p_0}(\xi_{p_0}),Q_{p_0}(Df_{p_0}^*(\tilde{\xi}_{f(p_0)})) \big ). 
\end{align*} 
Since $Q_{p_0}$ is an isometry, we have 
\[(\text{\rm Hess} \, h)_{p_0} \big ( Q_{p_0}^*(\tilde{\xi}_{f(p_0)}),Q_{p_0}^*(\tilde{\xi}_{f(p_0)}) \big ) 
\geq \theta \, |Q_{p_0}^*(\tilde{\xi}_{f(p_0)})|^2 = \theta \, |\tilde{\xi}_{f(p_0)}|^2 = \theta.\] 
Moreover, we have 
\[(\text{\rm Hess} \, \tilde{h})_{f(p_0)} \big ( Q_{p_0}(\xi_{p_0}),Q_{p_0}(Df_{p_0}^*(\tilde{\xi}_{f(p_0)})) \big ) \geq 0\] 
since $\tilde{h}$ is convex and $Df_{p_0}^*(\tilde{\xi}_{f(p_0)})$ is a positive multiple of $\xi_{p_0}$. Putting these facts together, we obtain 
\[\Big ( \frac{C_1}{\theta} - \mu \Big ) \, \varphi(p_0,f(p_0)) \geq \theta.\] 
Since $\varphi(p_0,f(p_0)) \geq 0$ and $\mu \geq 0$, we conclude that $\varphi(p_0,f(p_0)) \geq \frac{\theta^2}{C_1}$. On the other hand, we have $\inf_{p \in \partial \Omega} \varphi(p,f(p)) = \varphi(p_0,f(p_0))$ by definition of $p_0$. From this, the assertion follows. \\

\begin{proposition} 
\label{differential.of.f}
Suppose that $p$ is a point in $\Omega$ such that $\tilde{\xi}_{f(p)} \neq 0$, and $\{v_1,v_2\}$ is an orthonormal basis of $T_p N$. Let 
\[\Gamma(p) = \sum_{k=1}^2 \langle Df_p(v_k),\tilde{\xi}_{f(p)} \rangle \, \langle Q_p(v_k),\tilde{\xi}_{f(p)} \rangle > 0.\] 
Then we have 
\[\langle Df_p(v_1),Q_p(v_1) \rangle = \frac{1}{\Gamma(p)} \, \big ( \langle Df_p(v_1),\tilde{\xi}_{f(p)} \rangle^2 + \langle Q_p(v_2),\tilde{\xi}_{f(p)} \rangle^2 \big )\] 
and
\[\langle Df_p(v_2),Q_p(v_2) \rangle = \frac{1}{\Gamma(p)} \, \big ( \langle Df_p(v_2),\tilde{\xi}_{f(p)} \rangle^2 + \langle Q_p(v_1),\tilde{\xi}_{f(p)} \rangle^2 \big ).\] 
Moreover, we have 
\begin{align*} 
&\langle Df_p(v_1),Q_p(v_2) \rangle = \langle Df_p(v_2),Q_p(v_1) \rangle \\ 
&= \frac{1}{\Gamma(p)} \, \big ( \langle Df_p(v_1),\tilde{\xi}_{f(p)} \rangle \, \langle Df_p(v_2),\tilde{\xi}_{f(p)} \rangle - \langle Q_p(v_1),\tilde{\xi}_{f(p)} \rangle \, \langle Q_p(v_2),\tilde{\xi}_{f(p)} \rangle \big ). 
\end{align*}
\end{proposition}

\textbf{Proof.} 
By definition of $Q_p$, we have 
\[\langle Df_p(v_k),Q_p(v_l) \rangle = \big \langle \big [ Df_p^* \, Df_p \big ]^{\frac{1}{2}} v_k,v_l \big \rangle\]
for $1 \leq k,l \leq 2$. Hence, the matrix $\langle Df_p(v_k),Q_p(v_l) \rangle$, $1 \leq k,l \leq 2$, is positive definite with 
determinant $1$. Using the chain rule, we obtain 
\[\sum_{l=1}^2 \langle Df_p(v_k),Q_p(v_l) \rangle \, \langle Q_p(v_l),\tilde{\xi}_{f(p)} \rangle 
= \langle Df_p(v_k),\tilde{\xi}_{f(p)} \rangle\] 
for $k = 1,2$. The assertion follows now from a straightforward calculation. \\

\begin{corollary} 
\label{bound.for.df}
There exists a constant $C_3$, depending only on $h$ and $\tilde{h}$, such that 
\[\det(I + Df_p^* \, Df_p) \leq C_3\] 
for all points $p \in \partial \Omega$.
\end{corollary} 

\textbf{Proof.} 
Fix a point $p \in \partial \Omega$. Let $v_1$ be the outward-pointing unit normal vector to $\partial \Omega$ at $p$, and let $v_2$ be a unit vector tangential to $\partial \Omega$. Since $Df_p(v_2)$ is a tangent vector to $\partial \tilde{\Omega}$, we have $\langle Df_p(v_2),\tilde{\xi}_{f(p)} \rangle = 0$. Using Proposition \ref{differential.of.f}, we obtain 
\[\langle Df_p(v_1),Q_p(v_1) \rangle = \frac{\langle Df_p(v_1),\tilde{\xi}_{f(p)} \rangle^2 + \langle Q_p(v_2),\tilde{\xi}_{f(p)} \rangle^2}{\langle Df_p(v_1),\tilde{\xi}_{f(p)} \rangle \, \langle Q_p(v_1),\tilde{\xi}_{f(p)} \rangle},\] 
\[\langle Df_p(v_2),Q_p(v_2) \rangle = \frac{\langle Q_p(v_1),\tilde{\xi}_{f(p)} \rangle}{\langle Df_p(v_1),\tilde{\xi}_{f(p)} \rangle},\] 
and 
\[\langle Df_p(v_1),Q_p(v_2) \rangle = \langle Df_p(v_2),Q_p(v_1) \rangle 
= -\frac{\langle Q_p(v_2),\tilde{\xi}_{f(p)} \rangle}{\langle Df_p(v_1),\tilde{\xi}_{f(p)} \rangle}.\] 
We claim that $\langle Df_p(v_k),Q_p(v_l) \rangle$ is uniformly bounded from above for all $k,l$. By Corollary \ref{estimate.for.the.boundary.defining.functions.2}, we have 
\[\theta^2 \leq \langle Df_p(v_1),\tilde{\xi}_{f(p)} \rangle \leq \frac{1}{\theta^2}.\] 
Moreover, it follows from Proposition \ref{lower.bound.for.varphi} that 
\[\langle Q_p(v_1),\tilde{\xi}_{f(p)} \rangle = \varphi(p,f(p)) \geq \frac{\theta^2}{C_1}.\] 
Hence, there exists a constant $C_2$ such that 
\[|\langle Df_p(v_k),Q_p(v_l) \rangle| \leq C_2\] 
for all $k,l$. Thus, we conclude that 
\[\det(I + Df_p^* \, Df_p) = 2 + \sum_{k,l=1}^2 \langle Df_p(v_k),Q_p(v_l) \rangle^2 \leq 2 + 4C_2^2.\]

\section{The interior gradient estimate}

In this section, we show that the singular values of $Df_p$ are uniformly bounded for all $p \in \Omega$. To that end, we define a function $\beta: \Sigma \to \mathbb{R}$ by 
\[\beta(p,f(p)) = \frac{2}{\sqrt{\det(I + Df_p^* \, Df_p)}}.\]
It is easy to see that $0 < \beta(p,f(p)) \leq 1$ for all $p \in \Omega$. \\

\begin{proposition}
\label{global.estimate.for.df}
There exists a constant $C_4$, depending only on $h$ and $\tilde{h}$, such that 
\[\det(I + Df_p^* \, Df_p) \leq C_4\] 
for all points $p \in \Omega$.
\end{proposition}

\textbf{Proof.}
Since $\text{\rm Ric}_M = \kappa \, g_M$, it follows from work of M.T.~Wang that
\[\Delta_\Sigma \beta + 2\beta \, \sum_{k,l=1}^2 |I\!I(e_k,e_l)|^2 + \kappa \, \beta \, (1 - \beta^2) = 0,\] 
where $\{e_1,e_2\}$ is an orthonormal basis for $T\Sigma$.
(This follows from equation (3.9) in \cite{Wang2};
see also \cite{Wang1}, equation (2.2).) This implies
\[\Delta_\Sigma(\log \beta) + |\nabla^\Sigma(\log \beta)|^2 + 2 \, \sum_{k,l=1}^2 |I\!I(e_k,e_l)|^2 + \kappa \, (1 - \beta^2) = 0,\] 
hence 
\[\Delta_\Sigma(\log \beta) + \kappa \leq 0.\] 
On the other hand, we have $\Delta_\Sigma H \geq \theta > 0$ by Proposition \ref{laplacian.of.H}. Therefore, the function $\log \beta + \frac{\kappa}{\theta} \, H$ is superharmonic. Using the maximum principle, we obtain 
\begin{align*} 
&\sup_{p \in \Omega} \Big ( \log \det(I + Df_p^* \, Df_p) - \frac{2\kappa}{\theta} \, h(p) \Big ) \\ 
&= \sup_{p \in \partial \Omega} \log \det(I + Df_p^* \, Df_p) \leq \log C_3. 
\end{align*} 
Thus, we conclude that 
\[\det(I + Df_p^* \, Df_p) \leq C_3 \, \exp \bigg ( \frac{2\kappa}{\theta} \, \inf_\Omega h \bigg )\] 
for all $p \in \Omega$.

\section{Estimates in $C^{1,\alpha}$}

In this section, we prove uniform estimates for $f$ in $C^{1,\alpha}$. In a first step, we will establish uniform $C^{1,\alpha}$ bounds for $f$ in a neighborhood of $\partial \Omega$. 

\begin{lemma} 
\label{pde}
Assume that $F: \tilde{\Omega} \to \mathbb{R}$ is a smooth function. Then the function $F \circ f: \Omega \to \mathbb{R}$ satisfies 
\begin{align*} 
&\text{\rm tr} \Big [ (I + Df_p^* \, Df_p)^{-1} \, (\text{\rm Hess} \, (F \circ f))_p \Big ] \\ 
&= \text{\rm tr} \Big [ Df_p \, (I + Df_p^* \, Df_p)^{-1} \, Df_p^* \, (\text{\rm Hess} \, F)_{f(p)} \Big ] 
\end{align*} 
for all $p \in \Omega$.
\end{lemma}

\textbf{Proof.} 
We define a function $G: \Omega \times \tilde{\Omega} \to \mathbb{R}$ by $G(p,q) = F(f(p)) - F(q)$ for all points $(p,q) \in \Omega \times \tilde{\Omega}$. Clearly, $G|_\Sigma = 0$, hence $\Delta_\Sigma(G|_\Sigma) = 0$. Since $\Sigma$ is a minimal submanifold of $M$, we obtain 
\[\sum_{k=1}^2 (\text{\rm Hess}_M \, G)_{(p,f(p))}(e_k,e_k) = 0,\] 
where $\{e_1,e_2\}$ is an orthonormal basis for the tangent space $T_{(p,f(p))} \Sigma$. From this, the assertion follows easily. \\

\begin{proposition} 
\label{regularity.1}
There exist positive constants $\alpha$ and $C_5$ such that $[\tilde{h} \circ f]_{C^{1,\alpha}(\Omega)} \leq C_5$. 
\end{proposition}

\textbf{Proof.} 
It follows from the previous lemma that 
\[\theta \leq \text{\rm tr} \Big [ (I + Df_p^* \, Df_p)^{-1} \, (\text{\rm Hess} \, (\tilde{h} \circ f))_p \Big ] \leq \frac{1}{\theta}\] 
for all $p \in \Omega$. By Proposition \ref{global.estimate.for.df}, the eigenvalues 
of the symmetric operator $I + Df_p^* \, Df_p: T_p N \to T_p N$ are uniformly bounded from above and below. Since 
$\tilde{h} \circ f$ vanishes along the boundary of $\Omega$, the assertion follows from work of Morrey and Nirenberg (see \cite{Gilbarg-Trudinger}, Section 12.2, pp. 300--304). \\

In order to obtain uniform bounds for $f$ in $C^{1,\alpha}$, we choose a globally defined orthonormal frame $\{v_1,v_2\}$ on $N$. For abbreviation, let 
\begin{align*}
&a_{kl}(p) = \langle Q_p(v_k),(v_l)_{f(p)} \rangle \\ 
&b_k(p) = \langle Df_p(v_k),\tilde{\xi}_{f(p)} \rangle = \langle v_k,\nabla (\tilde{h} \circ f)_p \rangle \\ 
&c_l(p) = \langle (v_l)_{f(p)},\tilde{\xi}_{f(p)} \rangle 
\end{align*} 
for $1 \leq k,l \leq 2$. The following result implies that the gradient of $a_{kl}$ is uniformly bounded:

\begin{lemma} 
\label{q}
The gradient of the function $a_{kl}$ is given by 
\[\langle \nabla a_{kl},v_j \rangle = \langle Q_p(\nabla_{v_j} v_k),v_l \rangle + \langle Q_p(v_k),\nabla_{Df_p(v_j)} v_l \rangle\] 
for $j = 1,2$.
\end{lemma}

\textbf{Proof.} 
This follows from the same arguments that we used in the proof Proposition \ref{properties.of.varphi}. \\

It follows from Proposition \ref{global.estimate.for.df} that the eigenvalues of $Df_p^* \, Df_p$ lie in the interval $[\frac{1}{C_4},C_4]$ for all $p \in \Omega$. This implies 
\begin{align*}
\Gamma(p) 
&= \sum_{k=1}^2 \langle Df_p(v_k),\tilde{\xi}_{f(p)} \rangle \, \langle Q_p(v_k),\tilde{\xi}_{f(p)} \rangle \\ 
&= \big \langle Df_p^*(\tilde{\xi}_{f(p)}),Q_p^*(\tilde{\xi}_{f(p)}) \big \rangle \\
&= \big \langle [Df_p^* \, Df_p]^{\frac{1}{2}} \, Q_p^*(\tilde{\xi}_{f(p)}),Q_p^*(\tilde{\xi}_{f(p)}) \big \rangle \\ 
&\geq \frac{1}{\sqrt{C_4}} \, |Q_p^*(\tilde{\xi}_{f(p)})|^2 \\ 
&= \frac{1}{\sqrt{C_4}} \, |\tilde{\xi}_{f(p)}|^2 
\end{align*} 
for all $p \in \Omega$. 

By Proposition \ref{boundary.defining.functions}, we have $|\tilde{\xi}_q| = 1$ for all points $q \in \partial \tilde{\Omega}$. Hence, we can find a positive real number $\rho$ such that $|\tilde{\xi}_q| \geq \frac{1}{2}$ for all points $q \in \tilde{\Omega}$ satisfying $\tilde{h}(q) \geq -\rho$. Let 
\[\Omega_1 = \{p \in \Omega: h(p) \geq -\theta^2 \rho\}\] 
and 
\[\Omega_2 = \{p \in \Omega: h(p) \leq -\frac{1}{2} \, \theta^2 \rho\}.\]
By Proposition \ref{estimate.for.the.boundary.defining.functions.1}, we have $\tilde{h}(f(p)) \geq -\rho$ for all points $p \in \Omega_1$. This implies $|\tilde{\xi}_{f(p)}| \geq \frac{1}{2}$ for all $p \in \Omega_1$. Putting these facts together, we obtain 
\[\Gamma(p) \geq \frac{1}{\sqrt{C_4}} \, |\tilde{\xi}_{f(p)}|^2 \geq \frac{1}{4\sqrt{C_4}}\] 
for all points $p \in \Omega_1$.

\begin{proposition}
\label{regularity.2}
There exists a constant $C_9$ such that $[f]_{C^{1,\alpha}(\Omega_1)} \leq C_9$.
\end{proposition}

\textbf{Proof.}
Consider the functions $\chi_{kl}(p) = \langle Df_p(v_k),(v_l)_{f(p)} \rangle$ ($1 \leq k,l \leq 2$). Using Proposition \ref{differential.of.f} and the identity $a_{11} \, a_{22} - a_{12} \, a_{21} = 1$, we obtain 
\begin{align*}
\chi_{11} &= \frac{1}{\Gamma} \, \Big [ (a_{11} \, b_1 + a_{21} \, b_2) \, b_1 + (a_{21} \, c_1 + a_{22} \, c_2) \, c_2 \Big ] \\ 
\chi_{12} &= \frac{1}{\Gamma} \, \Big [ (a_{12} \, b_1 + a_{22} \, b_2) \, b_1 - (a_{21} \, c_1 + a_{22} \, c_2) \, c_1 \Big ] \\ 
\chi_{21} &= \frac{1}{\Gamma} \, \Big [ (a_{11} \, b_1 + a_{21} \, b_2) \, b_2 - (a_{11} \, c_1 + a_{12} \, c_2) \, c_2 \Big ] \\ 
\chi_{22} &= \frac{1}{\Gamma} \, \Big [ (a_{12} \, b_1 + a_{22} \, b_2) \, b_2 + (a_{11} \, c_1 + a_{12} \, c_2) \, c_1 \Big ]. 
\end{align*} 
Moreover, the function $\Gamma$ can be written in the form 
\[\Gamma = \sum_{k,l=1}^2 a_{kl} \, b_k \, c_l.\] 
Since $\Gamma(p) \geq \frac{1}{4\sqrt{C_4}}$ for all $p \in \Omega_1$, we can find a constant $C_6$ such that
\begin{align*} 
\sum_{k,l=1}^2 [\chi_{kl}]_{C^\alpha(\Omega_1)} 
&\leq C_6 \, \bigg ( \sum_{k,l=1}^2 [a_{kl}]_{C^\alpha(\Omega_1)} + \sum_{k=1}^2 [b_k]_{C^\alpha(\Omega_1)} + \sum_{l=1}^2 [c_l]_{C^\alpha(\Omega_1)} \bigg ). 
\end{align*} 
It follows from Lemma \ref{q} that the gradient of $a_{kl}$ is uniformly bounded for all $1 \leq k,l \leq 2$. Moreover, it is easy to see that the gradient of the function $c_l$ is uniformly bounded for $l = 1,2$. Finally, we have 
\[\sum_{k=1}^2 [b_k]_{C^\alpha(\Omega)} \leq C_7\] 
by Proposition \ref{regularity.1}. Putting these facts together, we obtain 
\[\sum_{k,l=1}^2 [\chi_{kl}]_{C^\alpha(\Omega_1)} \leq C_8\] 
for some constant $C_8$. From this, the assertion follows. \\

\begin{proposition}
\label{regularity.3} 
Assume that $\alpha$ is sufficiently small. Then there exists a constant $C_{12}$ such that 
$[f]_{C^{1,\alpha}(\Omega_2)} \leq C_{12}$. 
\end{proposition}

\textbf{Proof.} 
Fix a global coordinate system on $N$, and let $f_1,f_2: \Omega \to \mathbb{R}$ be the coordinate functions of $f$. Using Lemma \ref{pde}, we obtain 
\[\Big | \text{\rm tr} \Big [ (I + Df_p^* \, Df_p)^{-1} \, (\text{\rm Hess} \, f_j)_p \Big ] \Big | \leq C_{10}\] 
for all $p \in \Omega$ and $j = 1,2$. Hence, by Theorem 12.4 in \cite{Gilbarg-Trudinger}, we can find a constant $C_{11}$ such that $[f_j]_{C^{1,\alpha}(\Omega_2)} \leq C_{11}$ for $j = 1,2$. From this, the assertion follows. \\

It follows from the preceeding arguments that $f$ is bounded in $C^{1,\alpha}(\Omega)$. In order to prove higher regularity, we proceed as follows: assume that $f$ is bounded in $C^{m,\alpha}(\Omega)$ for some positive integer $m$. It follows from Lemma \ref{pde} and Schauder theory that the function $\tilde{h} \circ f$ is bounded in $C^{m+1,\alpha}(\Omega)$. Hence, the function $b_k$ is bounded in $C^{m,\alpha}(\Omega)$ for $k = 1,2$. Moreover, Lemma \ref{q} implies that the gradient of $a_{kl}$ is bounded in $C^{m-1,\alpha}(\Omega)$ for $1 \leq k,l \leq 2$. Therefore, the function $a_{kl}$ is bounded in $C^{m,\alpha}(\Omega)$ for $1 \leq k,l \leq 2$. Finally, it is easy to see that the function $c_l$ is bounded in $C^{m,\alpha}(\Omega)$ for $l = 1,2$. Consequently, the function $\chi_{kl}(p) = \langle Df_p(v_k),(v_l)_{f(p)} \rangle$ is bounded in $C^{m,\alpha}(\Omega_1)$ for $1 \leq k,l \leq 2$. On the other hand, it follows from interior Schauder estimates that $f$ is bounded in $C^{m+1,\alpha}(\Omega_2)$. Putting these facts together, we conclude that $f$ is bounded in $C^{m+1,\alpha}(\Omega)$.

\section{The linearized operator}

As above, we fix two strictly convex domains $\Omega,\tilde{\Omega} \subset N$ with smooth boundary. Moreover, we fix two points $\overline{p} \in \partial \Omega$ and $\overline{q} \in \partial \tilde{\Omega}$. Suppose that $f: \Omega \to \tilde{\Omega}$ is a minimal Lagrangian diffeomorphism satisfying $f(\overline{p}) = \overline{q}$. We claim that the linearized operator at $f$ is invertible. 

In order to prove this, we fix a real number $\alpha \in (0,1)$. We denote by $\mathcal{M}$ the space of all diffeomorphisms $\varphi: \Omega \to \Omega$ of class $C^{3,\alpha}$ that are area-preserving and orientation-preserving. It follows from the implicit function theorem that $\mathcal{M}$ is a Banach manifold. The tangent space to $\mathcal{M}$ at the identity can be identified with the space of all divergence-free vector fields on $\Omega$ of class $C^{3,\alpha}$ that are tangential at the boundary $\partial \Omega$. We will denote this space by $\mathcal{X}$. In other words, $\mathcal{X}$ consists of all vector fields of the form $J \, \nabla u$, where $u: \Omega \to \mathbb{R}$ is a function of class $C^{4,\alpha}$ satisfying $u|_{\partial \Omega} = 0$. Finally, we denote by $\mathcal{Y}$ the space of all closed one-forms on $\Omega$ of class $C^{1,\alpha}$.

We define a map $\mathcal{H}: \mathcal{M} \to \mathcal{Y}$ as follows: for each $\varphi \in \mathcal{M}$, we denote by $\mathcal{H}(\varphi)$ the mean curvature one-form associated with the Lagrangian embedding $p \mapsto (\varphi(p),f(p))$. Note that $\mathcal{H}(\varphi) \in \mathcal{Y}$ since the mean curvature one-form associated with a Lagrangian embedding is closed.

\begin{proposition}
\label{linearized.operator.1}
The linearized operator $D\mathcal{H}_{\text{\rm id}}: \mathcal{X} \to \mathcal{Y}$ sends a vector field $J \, \nabla u \in \mathcal{X}$ to the one-form $d(\Delta_g u + \kappa \, u) \in \mathcal{Y}$. 
Here, $g$ denotes the pull-back of the product metric on $M$ under the map $p \mapsto (p,f(p))$. 
\end{proposition}

\textbf{Proof.} 
Consider a one-parameter family of diffeomorphisms $\varphi_s \in \mathcal{M}$ such that $\varphi_0 = \text{\rm id}$ and 
$\frac{d}{ds} \varphi_s \big |_{s=0} = J \, \nabla u$. We define a one-parameter family of Lagrangian embeddings $F_s: \Omega \to M$ by \[F_s: p \mapsto (\varphi_s(p),f(p)).\] We denote by $V$ the variation vector field associated with this family of Lagrangian embeddings. This vector field is given by $V = \frac{d}{ds} F_s \big |_{s=0} = (J \, \nabla u,0)$. We next define a one-form $\eta$ on $\Omega$ by $\eta(w) = -\langle JV,DF_0(w) \rangle$. Since $-JV = (\nabla u,0)$ and $DF_0(w) = (w,Df(w))$, we have $\eta(w) = \langle \nabla u,w \rangle$, hence $\eta = du$. Using Proposition \ref{linearized.equation}, we obtain 
\[\frac{d}{ds} \mathcal{H}(\varphi_s) \big |_{s=0} = -d\delta_g \eta + \kappa \, \eta = d(-\delta_g du + \kappa \, u).\] 
Since $\Delta_g u = -\delta_g du$, the assertion follows. \\

We next define a map $\mathcal{G}: \mathcal{M} \to \mathcal{Y} \times \partial \Omega$ by 
$\mathcal{G}(\varphi) = \big ( \mathcal{H}(\varphi), \, \varphi(\overline{p}) \big )$. Note that $\mathcal{G}(\varphi) = (0,\overline{p})$ if and only if $f \circ \varphi^{-1}: \Omega \to \tilde{\Omega}$ is a minimal Lagrangian diffeomorphism that maps $\overline{p}$ to $\overline{q}$.

\begin{proposition} 
\label{linearized.operator.2}
The linearized operator $D\mathcal{G}_{\text{\rm id}}: \mathcal{X} \to \mathcal{Y} \times T_{\overline{p}}(\partial \Omega)$ is invertible.
\end{proposition}

\textbf{Proof.} 
The linearized operator $D\mathcal{G}_{\text{\rm id}}: \mathcal{X} \to \mathcal{Y} \times T_{\overline{p}}(\partial \Omega)$ sends a vector field $J \, \nabla u \in \mathcal{X}$ to the pair 
\[\big ( d(\Delta_g u + \kappa \, u), \, J \, \nabla u_{\overline{p}} \big ) \in \mathcal{Y} \times T_{\overline{p}} (\partial \Omega).\] 
Since $u|_{\partial \Omega} = 0$, the vector field $J \, \nabla u$ is tangential to the boundary $\partial \Omega$. 

We claim that the operator $D\mathcal{G}_{\text{\rm id}}: \mathcal{X} \to \mathcal{Y} \times T_{\overline{p}}(\partial \Omega)$ is one-to-one. To prove this, suppose that $u$ is a real-valued function of class $C^{1,\alpha}$ such that $d(\Delta_g u + \kappa \, u) = 0$ in $\Omega$, $u = 0$ on $\partial \Omega$, and $\nabla u = 0$ at $\overline{p}$. This implies $\Delta_g u + \kappa \, u = c$ for some constant $c \in \mathbb{R}$. If the constant $c$ is positive, then $u$ is strictly negative in the interior of $\Omega$ by the maximum principle. Hence, the Hopf boundary point lemma (cf. \cite{Gilbarg-Trudinger}, Lemma 3.4) implies that the outer normal derivative of $u$ at $\overline{p}$ is strictly positive. 
This contradicts the fact that $\nabla u = 0$ at $\overline{p}$. Thus, we conclude that $c \leq 0$. An analogous argument shows that $c \geq 0$. Consequently, we must have $c = 0$. Using the maximum principle, we deduce that $u = 0$. Thus, the operator $D\mathcal{G}_{\text{\rm id}}: \mathcal{X} \to \mathcal{Y} \times T_{\overline{p}}(\partial \Omega)$ is one-to-one.

A similar argument shows that $D\mathcal{G}_{\text{\rm id}}: \mathcal{X} \to \mathcal{Y} \times T_{\overline{p}}(\partial \Omega)$ is onto. This completes the proof. \\

\section{The continuity method}

In this section, we prove Theorem \ref{existence.uniqueness} using the continuity method. To that end, we deform $\Omega$ and $\tilde{\Omega}$ to the flat unit disk $\mathbb{B}^2 \subset \mathbb{R}^2$. There is a convenient way of performing this deformation, which we describe next.

Let $h$ and $\tilde{h}$ be the boundary defining functions constructed in Section 2. For each $t \in (0,1]$, we consider the sub-level sets of $h$ and $\tilde{h}$ with area $t^2 \, \text{\rm area}(\Omega) = t^2 \, \text{\rm area}(\tilde{\Omega})$. More previsely, we define two functions $A,\tilde{A}: (0,1] \to (-\infty,0]$ by 
\begin{align*} 
&\text{\rm area}(\{p \in \Omega: h(p) \leq A(t)\}) = t^2 \, \text{\rm area}(\Omega) \\ 
&\text{\rm area}(\{q \in \tilde{\Omega}: \tilde{h}(q) \leq \tilde{A}(t)\}) = t^2 \, \text{\rm area}(\tilde{\Omega}) 
\end{align*} 
for $t \in (0,1]$. For each $t \in (0,1]$, we consider the domains 
\begin{align*} 
\Omega_t &= \{p \in \Omega: h(p) \leq A(t)\} \\ 
\tilde{\Omega}_t &= \{q \in \tilde{\Omega}: \tilde{h}(q) \leq \tilde{A}(t)\}. 
\end{align*} 
It is easy to see that $\Omega_t$ and $\tilde{\Omega}_t$ are strictly convex domains with smooth boundary. Moreover, $\Omega_t$ and $\tilde{\Omega}_t$ have the same area. It follows from results in Section 2 that $\Omega_t$ and $\tilde{\Omega}_t$ are geodesic disks if $t > 0$ is sufficiently small. For each $t \in (0,1]$, we consider the following problem: \\

$(\star_t)$ \textit{Find all minimal Lagrangian diffeomorphisms $f: \Omega_t \to \tilde{\Omega}_t$ that map a given point on the boundary of $\Omega_t$ to a given point on the boundary of $\tilde{\Omega}_t$.} \\

We now pass to the limit as $t \to 0$. After suitable rescaling, the domains $\Omega_t$ and $\tilde{\Omega}_t$ converge to the flat unit disk $\mathbb{B}^2$. Hence, if we send $t \to 0$, then the problem $(\star_t)$ reduces to the following problem: \\

$(\star_0)$ \textit{Find all minimal Lagrangian diffeomorphisms $f: \mathbb{B}^2 \to \mathbb{B}^2$ that map one given point on $\partial \mathbb{B}^2$ to another given point on $\partial \mathbb{B}^2$.} \\

We claim that the problem $(\star_0)$ has a unique solution. To prove this, we need a uniqueness result for the second boundary value problem for the Monge-Amp\`ere equation:

\begin{proposition} 
\label{uniqueness.for.monge.ampere}
Suppose that $u,v: \mathbb{B}^2 \to \mathbb{R}$ are smooth convex functions satisfying $\det D^2 u(x) = \det D^2 v(x) = 1$ for all $x \in \mathbb{B}^2$. Moreover, suppose that the gradient mappings $x \mapsto \nabla u(x)$ and $x \mapsto \nabla v(x)$ map $\mathbb{B}^2$ to itself. Then the function $u(x) - v(x)$ is constant.
\end{proposition}

\textbf{Proof.} 
This is a subcase of a general uniqueness result due to Y.~Brenier \cite{Brenier}. A proof based on PDE methods was given by P.~Delano\"e \cite{Delanoe}. \\

\begin{proposition} 
\label{classification}
Suppose that $f$ is a minimal Lagrangian diffeomorphism from the flat unit disk $\mathbb{B}^2$ to itself. Then $f$ is a rotation.
\end{proposition}

\textbf{Proof.} 
Since the graph of $f$ is a minimal surface, there exists a constant $\gamma \in \mathbb{R}$ such that \[\cos \gamma \, (\partial_1 f_2(x) - \partial_2 f_1(x)) = \sin \gamma \, (\partial_1 f_1(x) + \partial_2 f_2(x))\] for all $x \in \mathbb{B}^2$. Hence, there exists a smooth function $u: \mathbb{B}^2 \to \mathbb{R}$ such that 
\begin{align*} 
\partial_1 u(x) &= \cos \gamma \, f_1(x) + \sin \gamma \, f_2(x) \\ 
\partial_2 u(x) &= -\sin \gamma \, f_1(x) + \cos \gamma \, f_2(x). 
\end{align*} 
By assumption, $f$ is a diffeomorphism from $\mathbb{B}^2$ to itself. Hence, the gradient mapping $x \mapsto \nabla u(x)$ is a diffeomorphism from $\mathbb{B}^2$ to itself. Since $f$ is area-preserving and orientation-preserving, we have $\det D^2 u(x) = \det Df(x) = 1$ for all $x \in \mathbb{B}^2$. Consequently, the function $u$ is either convex or concave. If $u$ is convex, it follows from Proposition \ref{uniqueness.for.monge.ampere} that $\frac{1}{2} \, |x|^2 - u(x)$ is constant. This implies 
\begin{align*} 
f_1(x) &= \cos \gamma \, \partial_1 u(x) - \sin \gamma \, \partial_2 u(x) = \cos \gamma \, x_1 - \sin \gamma \, x_2 \\ 
f_2(x) &= \sin \gamma \, \partial_1 u(x) + \cos \gamma \, \partial_2 u(x) = \sin \gamma \, x_1 + \cos \gamma \, x_2. 
\end{align*} 
Similarly, if $u$ is concave, then Proposition \ref{uniqueness.for.monge.ampere} implies that $\frac{1}{2} \, |x|^2 + u(x)$ is constant. In this case, we obtain 
\begin{align*} 
f_1(x) &= \cos \gamma \, \partial_1 u(x) - \sin \gamma \, \partial_2 u(x) = -\cos \gamma \, x_1 + \sin \gamma \, x_2 \\ 
f_2(x) &= \sin \gamma \, \partial_1 u(x) + \cos \gamma \, \partial_2 u(x) = -\sin \gamma \, x_1 - \cos \gamma \, x_2. 
\end{align*} 
In either case, we conclude that $f$ is a rotation. \\

\begin{proposition}
For each $t \in (0,1]$, the problem $(\star_t)$ has exactly one solution.
\end{proposition}

\textbf{Proof.} By Theorem \ref{a.priori.estimates}, every minimal Lagrangian diffeomorphism from $\Omega_t$ to $\tilde{\Omega}_t$ is uniformly bounded in $C^m$ after rescaling. Moreover, it follows from Proposition \ref{linearized.operator.2} that every solution of $(\star_t)$ is non-degenerate. Consequently, $(\star_t)$ and $(\star_0)$ have the same number of solutions for all $t \in (0,1]$. Since $(\star_0)$ has a unique solution by Proposition \ref{classification}, the proof is complete.

\appendix

\section{The linearization of the Lagrangian minimal surface equation}

Let $M$ be a K\"ahler-Einstein manifold with $\text{\rm Ric}_M = \kappa \, g_M$. Moreover, let $F_s: \Sigma \to M$ be a one-parameter family of Lagrangian embeddings into $M$. For each $s$, we denote by $\mu_s$ the mean curvature one-form associated with $F_s$. Clearly, $\mu_s$ is a closed one-form on $\Sigma$. Finally, we define a one-form $\eta$ on $\Sigma$ by $\eta(X) = -\langle JV,DF_0(X) \rangle$ for $X \in T\Sigma$, where $V = \frac{\partial}{\partial s} F_s \big |_{s=0}$ denotes the variation vector field. 

\begin{proposition}
\label{linearized.equation}
Suppose that $F_0: \Sigma \to M$ is a minimal Lagrangian embedding, i.e. $\mu_0 = 0$. Then 
\[\frac{d}{ds} \mu_s \big |_{s=0} = -d\delta_g \eta + \kappa \, \eta.\] 
Here, $g$ denotes the pull-back of the Riemannian metric on $M$ under $F_0$.
\end{proposition}

\textbf{Proof.} 
Without loss of generality, we may assume that $\Sigma$ is a submanifold of $M$ and $F_0(p) = p$ for all $p \in \Sigma$. 
We can find a vector field $W \in T\Sigma$ such that $\eta(X) = -\langle JV,X \rangle = \langle W,X \rangle$ for all 
$X \in T\Sigma$. This implies $V - JW \in T\Sigma$, i.e. $JW$ is the normal component of $V$. The change of the mean curvature vector is given by the formula 
\begin{align*} 
&\sum_k \nabla_{e_k}^\perp \nabla_{e_k}^\perp (JW) 
+ \sum_k \big [ R_M(e_k,JW)e_k \big ]^\perp \\ 
&+ \sum_{k,l} \langle I\!I(e_k,e_l),JW \rangle \, I\!I(e_k,e_l) 
\end{align*}
(cf. \cite{Li}, Section 1, where a different sign convention for the curvature tensor is used). Hence, the change of the mean curvature one-form is given by 
\begin{align*} 
\frac{d}{ds} \mu_s(X) \big |_{s=0} 
&= \sum_k \langle \nabla_{e_k}^\perp \nabla_{e_k}^\perp (JW),JX \rangle + \sum_k R_M(e_k,JW,e_k,JX) \\ 
&+ \sum_{k,l} \langle I\!I(e_k,e_l),JW \rangle \, \langle I\!I(e_k,e_l),JX \rangle. 
\end{align*} 
for all $X \in T\Sigma$. Since $J$ is parallel, we have 
\[\nabla_{e_k}^\perp \nabla_{e_k}^\perp (JW) = J(\nabla_{e_k}^\Sigma \nabla_{e_k}^\Sigma W).\] 
Moreover, using the relation $\sum_k I\!I(e_k,e_k) = 0$ and the Gauss equations, we obtain 
\begin{align*} 
&\sum_k R_M(e_k,JW,e_k,JX) + \sum_{k,l} \langle I\!I(e_k,e_l),JW \rangle \, \langle I\!I(e_k,e_l),JX \rangle \\ 
&= \sum_k R_M(e_k,JW,e_k,JX) + \sum_{k,l} \langle I\!I(e_k,W),Je_l \rangle \, \langle I\!I(e_k,X),Je_l \rangle \\ 
&= \sum_k R_M(e_k,JW,e_k,JX) + \sum_k \langle I\!I(e_k,W),I\!I(e_k,X) \rangle \\ 
&= \sum_k R_M(e_k,JW,e_k,JX) + \sum_k R_M(e_k,W,e_k,X) - \sum_k R_\Sigma(e_k,W,e_k,X) \\ 
&= \sum_k R_M(Je_k,W,Je_k,X) + \sum_k R_M(e_k,W,e_k,X) - \sum_k R_\Sigma(e_k,W,e_k,X) \\ 
&= \text{\rm Ric}_M(W,X) - \text{\rm Ric}_\Sigma(W,X) \\
&= \kappa \, \langle W,X \rangle - \text{\rm Ric}_\Sigma(W,X) 
\end{align*}
for all $X \in T\Sigma$. Putting these facts together, we obtain 
\[\frac{d}{ds} \mu_s(X) \big |_{s=0} 
= \sum_k \langle \nabla_{e_k}^\Sigma \nabla_{e_k}^\Sigma W,X \rangle - \text{\rm Ric}_\Sigma(W,X) 
+ \kappa \, \langle W,X \rangle,\] 
hence 
\[\frac{d}{ds} \mu_s(X) \big |_{s=0} = \sum_k (\nabla_{e_k}^\Sigma \nabla_{e_k}^\Sigma \eta)(X) - \sum_k \text{\rm Ric}_\Sigma(e_k,X) \, \eta(e_k) + \kappa \, \eta(X)\] 
for $X \in T\Sigma$. On the other hand, we have 
\[(d\delta \eta)(X) + (\delta d\eta)(X) = -\sum_k (\nabla_{e_k}^\Sigma \nabla_{e_k}^\Sigma \eta)(X) + \sum_k \text{\rm Ric}_\Sigma(e_k,X) \, \eta(e_k)\] 
by the standard Bochner formula. Thus, we conclude that 
\[\frac{d}{ds} \mu_s \big |_{s=0} = -d\delta \eta - \delta d\eta + \kappa \, \eta.\] 
Finally, we have $d\eta = 0$ since $F_s$ is a one-parameter family of Lagrangian embeddings. From this, the assertion follows.

\end{document}